\documentclass{svjour3}
\usepackage{amsmath,amssymb,multirow}
\usepackage{bbm}
\usepackage{fullpage}
\usepackage{exscale,subfigure,booktabs}
\usepackage{color}
\usepackage{hyperref}
\usepackage{graphicx}
\usepackage[ruled,vlined]{algorithm2e}
\hypersetup{
    breaklinks=true,
	colorlinks=true,
	linkcolor=blue, % Couleur des liens internes
	citecolor=red, % Couleur des numŽros de la biblio dans le corps
	urlcolor=blue  } % Couleur des url
\usepackage[weather,alpine,misc,geometry]{ifsym}

\everymath{\displaystyle}
\newcommand{\la}{\lambda}

\newcommand {\R} {\mathbb R}
\newcommand {\N} {\mathbb N}

\newcommand{\norm}[1]{\left\Vert#1\right\Vert}
\newcommand{\abs}[1]{\left\vert#1\right\vert}
\newcommand{\ang}[1]{\left\langle #1 \right\rangle}

\newcounter{mycount}


\renewcommand {\theenumi} {\rm(\roman{enumi})}

\newcommand{\x}[1]{}
\usepackage{exscale,subfigure}
\usepackage{tcolorbox} %----------------------------------------------------------------
\usepackage{tikz}
\usepackage{pgfplots}
\usepgfplotslibrary{external}
\pgfplotsset{width=3cm,compat=1.9}
\newcommand{\ceil}[1]{\left \lceil{#1}\right \rceil }

\usepackage[natbibapa]{apacite}
\usepackage{breakcites}

\title{An exact cutting plane method for solving $p$-dispersion-sum problems}
\author{Sandy Spiers, Hoa T. Bui, Ryan Loxton}
\institute{
	Sandy Spiers (\Letter\,), Hoa T. Bui, Ryan Loxton
	\at
	ARC Training Centre for Transforming Maintenance through Data Science, and 
	Curtin Centre for Optimisation and Decision Science, Curtin University, Australia\\
	\email{sandy.spiers@postgrad.curtin.edu.au, hoa.bui@curtin.edu.au, r.loxton@curtin.edu.au}  
}
\date{Received: date / Accepted: date}
\usepackage{tikz}
\usetikzlibrary{patterns}
\usetikzlibrary{decorations.pathmorphing} % for snake lines
\usetikzlibrary{matrix} % for block alignment
\usetikzlibrary{arrows} % for arrow heads
\usetikzlibrary{calc} % for manimulation of coordinates
\usetikzlibrary{positioning,fit,calc} % used for the efficient
\usepackage{pgfplots}
\usepackage{pst-plot}
\begin{document}
	\maketitle
	\abstract{
		This paper aims to answer an open question recently posed in the literature, that is to find a fast exact method for solving the $p$-dispersion-sum problem (PDSP), a nonconcave quadratic binary maximization problem.
		We show that, since the Euclidean distance matrix defining the quadratic term in (PDSP) is always conditionally negative definite, the cutting plane method is exact for (PDSP) even in the absence of concavity.
		As such, the cutting plane method, which is primarily designed for concave maximisation problems, converges to the optimal solution of the (PDSP).
		%We establish that the cutting plane method can converge to the optimal solution of the (PDSP) without the need for convexification procedures commonly found when solving nonconcave quadratic problems with cutting planes.
		%We establish that without any convexification step, the cutting plane method applied to (PDSP) still converges to an optimal solution.
		%The $p$-dispersion problem is the problem of locating $p$ facilities at some of $n$ given locations so that the distance sum between the $p$ facilities is maximized. 
		The numerical results show that the method outperforms other exact methods for solving (PDSP), and can solve to optimality large instances of up to two thousand variables.}
	\keywords{ 
		Combinatorial optimization \and $p$-dispersion problem \and cutting planes \and branch and cut}
	
	\subclass{90-08; 90C26; 90C09; 90B80
		}
%	\setcounter{tocdepth}{4}
%	\setcounter{secnumdepth}{4}
%	\tableofcontents
	
	\section{Introduction}\label{intro}
	The problem of maximizing diversity arises in many practical settings. It involves selecting a subset of elements from a larger set to maximize some distance or dispersion metric. Since the conception of the maximum diversity problem in the 1980s by \cite{Kuby1987}, the interpretation of maximum diversity has taken many practical and theoretical forms. The topic has now reached a level of maturity where a multitude of problem variations, solution algorithms, and practical applications exist. Over the last thirty years, a significant quantity of research has focused on the max-sum problem \citep{Kuby1987}, which is to maximize the sum of distances between selected elements, and the max-min problem \citep{erkut_discrete_1990}, which is to maximize the minimum distance among selected points. For this paper, we focus our attention on the max-sum problem, which is referred to as the $p$-dispersion-sum problem (PDSP).
	
	Given a set of $n$ predefined locations, $v_1,\dots,v_n$ in a vector space $\mathbb{R}^s$ ($s\ge1$), the (PDSP) aims to find a subset of $p$ locations such that the sum of the distances between the $p$ points is maximized. Here, we consider $q_{ij}$ to be the distance between locations $i$ and $j$ defined by $q_{ij} = \left\| v_i - v_j \right\|^r$ where $r \in (0,2]$, and $\norm{\cdot}$ is the standard Euclidean distance in $\R^s$.  Note that when $r=1$, the distance $q_{ij}$ becomes the usual Euclidean distance between location $v_i$ and $v_j$.  Then, let $Q=[q_{ij}]$ denote the full %Euclidean 
	distance matrix where $i=1,\dots,n$ and $j=1\dots,n$. The (PDSP) is then given as
	\begin{alignat}{3}
		\tag{PDSP}
		\label{prob:pdsp}
		\max \quad & f(x) = \frac{1}{2}\ang{Qx,x}, && \\
		\label{cts:p}
		\text{s.t.} \quad
		& \sum_{i=1}^n x_i = p, &&\\
		\nonumber
		& x_i \in \{0,1\}, &\quad& 1 \le i \le n.
	\end{alignat}
	The \eqref{prob:pdsp} is known to be strongly NP-hard (see \cite{Kuo1993}, \cite{ravi_heuristic_1994}).  Even in the case where $r=2$ (square of the Euclidean distance), \cite{eremeev_maximum_2019} proved that the problem is still strongly NP-hard.
	
	The practical applications of the maximum dispersion problem are vast. One of the first examples presented in the literature is locating unwanted facilities on a network \citep{church_locating_1978}. Since then, many other researchers have relaxed the notion of distance to more general settings. One example is in genetics, where breeders attempt to maximize the diversity of traits among a breeding stock \citep{porter1975cowpea}. Furthermore, social diversity such as gender, cultural and ethnic diversity has become highly desirable in many communities, especially in a workplace setting \citep{roberge_recognizing_2010}. More recently, maximum dispersion has been used to find the optimal locations of chairs for COVID-19 social distancing \citep{ferrero-guillen_optimal_2022}.
	
	While research into heuristic and meta-heuristic approaches to the $p$-dispersion-sum problem has gathered significant interest (see \cite{Marti2022} for a recent review), the development of exact algorithms has fallen behind. One of the first exact approaches was presented in \cite{Kuo1993} and used linear reformulation techniques to transform the problem into an integer linear form. This was done in two ways. The first used a linearization technique presented in \cite{glover_technical_1974}, whereby the quadratic $x_i x_j$ terms are replaced by a new auxiliary variable $y_{ij}$.  The linear formulation of \eqref{prob:pdsp} is then given as
	\begin{alignat}{3}
		\label{prob:f2}
		\max \quad & \sum_{i=1}^{n-1} \sum_{j=i+1}^n q_{ij} y_{ij}, && \\
		\nonumber
		\text{s.t.} \quad
		& \sum_{i=1}^n x_i = p, &&\\
		\label{cts:y_1}
		& x_i + x_j - y_{ij} \le 1, &\quad& 1 \le i < j \le n,\\
		& -x_i + y_{ij} \le 0, &\quad& 1 \le i < j \le n,\\
		\label{cts:y_3}
		& -x_j + y_{ij} \le 0, &\quad& 1 \le i < j \le n,\\
		\nonumber
		& y_{ij} \ge 0, &\quad& 1 \le i < j \le n,\\
		\nonumber
		& x_i \in \{0,1\}, &\quad& 1 \le i \le n,
	\end{alignat}
	where constraints \eqref{cts:y_1}-\eqref{cts:y_3} enforce $y_{ij} = x_i x_j$.  
	A second reformulation that uses inequalities and real variables to handle quadratic terms, a technique first outlined in \cite{glover_improved_1975}, is given as
	\begin{alignat}{3}
		\label{prob:f3}
		\max \quad & \sum_{i=1}^{n-1} w_{i}, && \\
		\nonumber
		\text{s.t.} \quad
		& \sum_{i=1}^n x_i = p, &&\\
		\nonumber
		& w_{i}- x_i   \sum_{j=i+1}^n q_{ij}   \le 0, &\quad& 1 \le i \le n-1,\\
		\nonumber
		&  w_{i} - \sum_{j=i+1}^n q_{ij}  x_j \le 0, &\quad& 1 \le i \le n-1,\\
		\nonumber
		& w_{i} \ge 0, &\quad& 1 \le i \le n-1,\\
		\nonumber
		& x_i \in \{0,1\}, &\quad& 1 \le i \le n.
	\end{alignat}
	This formulation was shown in \cite{Marti2010} to be far more efficient than \eqref{prob:f2}. It was later used as the exact solver for the comprehensive empirical analyses presented in \cite{Parreno2021} and \cite{Marti2022}.
	
	The first significant advancement in exact methods for the \eqref{prob:pdsp} came in \cite{Pisinger2006}. The paper presented several upper bounds based on Lagrangean relaxation, semidefinite programming and reformulation techniques. The upper bounds are computationally cheap and can therefore be implemented in a branch and bound procedure. Numerical results show that for Euclidean distance problems, the procedure is capable of solving problems with $n=80$ with an average solve time of 60 seconds, but it struggles for sizes $n\ge100$. \cite{Marti2010} presented a branch and bound algorithm based on partial solutions, where a partial solution is a set of $k$ elements where $k < p$. Upper bounds are then calculated based on all other solutions that contain these $k$ elements. The objective function is split into three parts, and an upper bound for each is calculated. These bounds are then integrated into a branch and bound search tree. While the algorithm is faster than the linear formulation \eqref{prob:f3}, the numerical results show that it struggles to solve instances of $n=150$ in under an hour of computation time.
	
	This paper answers an open question posed in the recent review paper \cite{Marti2022}. That is, while progress in exact methods for the max-min diversity problem has advanced significantly, a fast exact solver for the $p$-dispersion-sum problem remains elusive. The max-sum problem remains the most widely studied problem variation, yet very few exact methods exist. One of the reasons for this might be that the problem is generally nonconcave, meaning the naive application of concave nonlinear programming techniques is not appropriate. However, when the distance measurements are taken as Euclidean, the problem exhibits certain special characteristics that allow for nonlinear programming techniques, particularly \emph{cutting plane methods}, to be applied, even in the presence of nonconcavity.

%n this paper, we apply the cutting plane method and its variants to solve the \eqref{prob:pdsp}. 
The cutting plane method (or outer approximation) (see \cite{duran1986outer}, \cite{leyffer1993deterministic}, \cite{yuan1988methode}) is one of several deterministic methods that provides general frameworks to tackle concave mixed integer problems. 
These methods require a concave objective function to guarantee convergence to optimality.
For nonconcave quadratic problems, the cutting plane algorithm requires an extra concave reformulation step before applying the cutting planes procedure. In particular, using the property that $x_i = x^2_i$ for $x_i\in \{0,1\}$, the nonconcave objective $f(x) = \tfrac{1}{2} \ang{Qx,x}$ is replaced by a concave function $f^\prime (x) := \tfrac{1}{2}\left( \ang{(Q-\la I_n)x,x}+\la\sum_{i=1}^nx_i \right)$, where $\la$ is the largest eigenvalue of $Q$ (see \cite{lima2017solution}), where here $I_n$ is the identity matrix of dimension $n$. This method has been implemented in commercial solvers
like \texttt{CPLEX} and \texttt{Gurobi} \citep{bliek1u2014solving,lima2017solution}. However, this approach is slow to converge (see \cite{bliek1u2014solving}, \cite{bonami2022classifier}). Hence, it is desirable to improve the efficiency of the cutting planes. We show in Section~\ref{basic} that  without the reformulation step the cutting plane algorithm applied to \eqref{prob:pdsp}, an inherently nonconcave problem, still converges to optimality.

The performance of the cutting plane algorithm is evaluated using the test set \texttt{GDK-d}, recently published in \cite{Parreno2021}. The authors of the recent review paper \cite{Marti2022} claim that this test set should replace the previously used Euclidean test sets as it considers very large test instances ($n=2000$) for the first time. Numerical results show that on the test set \texttt{GDK-d}, the cutting plane algorithm is vastly superior to other exact solvers. It can solve instances of up to $n=2000$ in an average of 10 seconds.  In many cases, the cutting plane algorithm is able to improve upon the best known heuristic and meta-heuristic solvers. However, a major shortcoming of the algorithm is its inability to solve instances generated using a large number of coordinates (i.e., large values of $s$).

The paper is organized as follows. In Section~\ref{basic}, we present an exact cutting plane approach for solving \eqref{prob:pdsp}. The convergence to optimality is established in Theorem~\ref{converge}. We then provide in Theorem~\ref{element}  an estimation of how many non-optimal solutions each cutting plane eliminates at each iteration.
	 Finally, in Section~\ref{NE}, we evaluate the effectiveness of the proposed cutting plane algorithm with the test set on the recently published test library.
	 
	\section{Cutting Plane Methodology}\label{basic}
	Denote  the feasible set of \eqref{prob:pdsp} as
	$$
	K := \left\{x\in\{0,1\}^n:\; \sum_{i=1}^nx_i = p \right\}.
	$$
	Let $h: \R^n\times \R^n \to \R$ be the tangent plane of $f$ defined as follows:
	\begin{equation}\label{h}
		h(x,y) :=\ang{\nabla f(y),x-y}+f(y),\quad \forall x\in \R^n,\quad\forall y\in \R^n.
	\end{equation}
	Given a set $A\subset K$, let $\Gamma_A \subset K\times \R$ be defined as
	$$
	\Gamma_A:=\left\{(x,\theta)\in \R^{n+1}:\; x\in K;\; \theta \le h(x,y),\; \forall y\in A \right\}.
	$$
	We consider the following auxiliary linear maximization problem
	\begin{equation}\label{LP_A}\tag{$\text{LP}_A$}
		\max_{(x,\theta)\in \Gamma_A} \theta.
	\end{equation}
	This linearization problem is known as the cutting-plane model of \eqref{prob:pdsp}, and can be written explicitly as 
	\begin{align}
	\notag	\max \quad &\theta\\
	\label{Ax-c1}	\text{s.t.}\quad & \theta \le h(x,y), \quad \forall y\in A,\\
	\notag	& x\in K.
	\end{align}
	Since it is essential to our analysis, before presenting the key results, we briefly discuss the properties of Euclidean distance matrices.
	An $n\times n$ matrix $D=[d_{ij}]$ ($n\ge 1$) is called an (squared) Euclidean distance matrix if there are $n$ vectors $v_1,\ldots,v_n$ in an Euclidean space $\R^s$ ($s\ge 1$) such that  $d_{ij} = \norm{v_i-v_j}^2$ for all $i,j=1,\ldots,n$, where $\norm{\cdot}$ is the Euclidean norm (see \cite{hayden1999methods}, \cite{schoenberg1937certain}, \cite{gower1982euclidean}). \cite{Schoenberg35remarksto,schoenberg1938metric} proved that a squared Euclidean distance matrix $D$ is \emph{conditionally negative definite}, that is $\ang{Dx,x} \le 0$ for any $x\in \R^n$ with $\sum_{i=1}^n x_i = 0$. Conversely, a symmetric nonnegative  matrix $D=[d_{ij}]$ with zero diagonal, i.e., $d_{ii} =0$, is an Euclidean distance matrix if and only if $D$ is conditionally negative definite (see {Theorem 4.7} of \cite{bapat1997nonnegative} and also \cite{hayden1999methods}, \cite{gower1982euclidean}). It is also proved in \cite{schoenberg1937certain}, and reproved in \cite{micchelli1984interpolation} (see also \cite{baxter1991conditionally}), that the $h$-exponential matrix $D^h$ of a Euclidean distance matrix $D$, defined as $D^h:= [d_{ij}^h]$ with $d_{ij}^h = \norm{v_i-v_j}^{2h}$, is also an Euclidean matrix if $h\in (0,1]$. In particular, when $h\le 1$, we can construct $n$ points $u_1,\ldots,u_n$ in the Euclidean space $\R^n$ such that 
	$$
	d_{ij}^h=\norm{v_i-v_j}^{2h} = \norm{u_i-u_j}^2,\quad\forall i,j=1,\ldots,n.
	$$
	The matrix $Q$ in Problem \eqref{prob:pdsp} is an Euclidean distance matrix and therefore it is conditionally negative definite.
	We now show that when $A = K$, the linear problem $(\text{LP}_K)$ is equivalent to the quadratic problem \eqref{prob:pdsp}.
	
	\begin{proposition}\label{P1}
	It holds that
	$$
	\max_{(x,\theta)\in \Gamma_K}\theta = \max_{x\in K} f(x).
	$$
	Furthermore, if $(x^*,\theta^*)$ is a solution of $(\text{LP}_K)$, then $x^*$ is a solution of \eqref{prob:pdsp}.
	\end{proposition}
\begin{proof}
	For any feasible solution $(x,\theta)$ of $(\text{LP}_K)$, we have $x\in K$ and
	$$
	\theta \le h(x,x) = f(x) \le \max_{z\in K} f(z).
	$$
	Therefore, $\max_{(x,\theta)\in \Gamma_K}\theta \le \max_{x\in K} f(x)$. Now, we prove the reverse inequality.
Taking into account that $Q$ is conditionally negative definite, and
for any feasible solutions $x,y\in K$ we have $\sum_{i=1}^n(x_i-y_i) =0$, the following inequality holds
$$\ang{Q(x-y),x-y} \le 0.$$
The inequality above yields
\begin{align*}
	h(x,y) - {f}(x) &= \ang{Qy,x-y}+ \tfrac{1}{2}\ang{Qy,y}- \tfrac{1}{2}\ang{Qx,x} \\
	& = \ang{Qy,x-y}+ \tfrac{1}{2} \ang{Q(y+x),y-x}\\
	& = -\tfrac{1}{2}\ang{Q(x-y),x-y} \ge 0.
\end{align*}
Thus, 
$$
h(x,y) \ge {f}(x),\quad \forall x,y\in K.
$$
Taking minimum over all $y\in K$ in the inequality above yields $\min_{y\in K}h(x,y) \ge f(x)$ for each $x\in K$.
Thus, $(x,\tilde{\theta}_x)$, with $\tilde{\theta}_x=\min_{y\in K}h(x,y)$, is a feasible solution of $(\text{LP}_K)$ such that $\tilde{\theta}_x\ge f(x)$, and therefore
	$$
\max_{(x,\theta)\in \Gamma_K}\theta \ge \max_{x\in K}\tilde{\theta}_x \ge \max_{x\in K} f(x).
$$
We have proved that $\max_{(x,\theta)\in \Gamma_K}\theta = \max_{x\in K} f(x)$.

Let $(x^*,\theta^*)$ be a solution of $(\text{LP}_K)$. 
Then, $x^*$ is feasible for \eqref{prob:pdsp}, and
$$
\max_{x\in K} f(x) = \theta^* \le h(x^*,x^*) = f(x^*),
$$
proving the second assertion.
\qed
\end{proof}

We propose a cutting plane approach to solve the quadratic problem \eqref{prob:pdsp} by solving the linear problem $(\text{LP}_K)$. Since it is not feasible to generate a cutting plane $h(x,y)\ge \theta$ for every $y\in K$, our cutting plane algorithm successively generates cuts of type \eqref{Ax-c1} whenever a candidate solution is found. Let $A_k$ denote the set of $A$ at each iteration $k$, and let $LB_k$ denote the lower bound at iteration $k$. We say $x\in K$ is a \emph{candidate solution} if there is $\theta > \text{LB}_k$ such that $(x,\theta)\in \Gamma_{A_k}$. Algorithm~\ref{Al-1} successively generates candidate solutions and adds the cutting planes to the linear model $(\text{LP}_{A_k})$ to eliminate non-optimal solutions until no candidate solution remains in the search space.

\begin{algorithm}[H]
	\SetAlgoLined
	\SetKwFor{Repeat}{repeat}{}{end}
	\vspace{0.5em}
	Take $x^0\in K$\\
	Set $A_0 \leftarrow \{x^0\}$,
	$\text{LB}_0$ $\leftarrow f(x_0)$, $k \leftarrow 1$ \\
	\While{$\exists(x^k,\theta^k) \in \Gamma_{A_{k-1}}$ s.t. $\theta^k > \text{\rm LB}_{k-1}$}{
	    $\text{LB}_k \leftarrow \max\{\text{LB}_{k-1},f(x^{k})\}$ \\
		$A_{k} \leftarrow A_{k-1}\cup\{x^k\}$ \\
	    $k \leftarrow k+1$\\
	}
	\label{Al-1}\caption{Cutting plane method for solving \eqref{prob:pdsp}.}
\end{algorithm}

~

This algorithm does not require solving the linear problem $(\text{LP}_{A_k})$ to optimality whenever an additional cut is added.  Rather, it looks for a feasible solution $(x^k,\theta^k) \in \Gamma_{A_{k-1}}$ that improves upon the current lower bound, i.e., $\theta^k > \text{\rm LB}_{k-1}$.   The implementation of Algorithm \ref{Al-1} is described in detail at the end of the section.  We now prove that Algorithm~\ref{Al-1} converges to an optimal solution of Problem \eqref{prob:pdsp}.

\begin{theorem}\label{converge}
	The sequence $\{x^k\}\subset K$ generated by Algorithm~\ref{Al-1} converges to an optimal solution of \eqref{prob:pdsp} after a finite number of steps.
\end{theorem}
\begin{proof}
	Consider the sequence  $\{x^k\}$ generated by Algorithm~\ref{Al-1}. We first prove that the algorithm terminates at an optimal solution.
	%We have $A_k\subset K$, and hence $\Gamma_{K} \subset \Gamma_{A_{k}}$, and the upper bound obtained during iteration $k$ is
	%	$$
	%	\text{UB}_k = \max_{(x,\theta)\in \Gamma_{A_{k}}} \theta \ge \max_{(x,\theta)\in \Gamma_{K} } \theta = \max_{x\in K} f(x)\ge \text{LB}_k.
	%	$$
	Suppose the algorithm terminates at step $k$, then for every $(\tilde{x},\tilde{\theta})\in \Gamma_{A_{k-1}}$, it holds that   
	%If $\text{UB}_k = \text{LB}_k$, then for every $(\tilde{x},\tilde{\theta})\in \Gamma_{A_{k+1}}$, from Proposition~\ref{P1}, it holds that     
	\begin{equation*}
	    \tilde{\theta} \le \text{LB}_{k-1} \le \max_{x\in K} f(x) = \max_{(x,\theta)\in \Gamma_{K}}\theta,
	\end{equation*}
	where the last equality follows Proposition~\ref{P1}.
	Taking the maximum over all $(\tilde{x},\tilde{\theta})\in \Gamma_{A_{k-1}}$ in the first inequality, we obtain
	\begin{equation}
	    \label{T1P1}
	    \max_{(x,\theta)\in \Gamma_{A_{k-1}}}\theta \le \text{LB}_{k-1}  \le \max_{(x,\theta)\in \Gamma_{K}}\theta
	\end{equation}
	Note that because $A_{k-1}\subset K$, we have $\Gamma_{K} \subset \Gamma_{A_{k-1}}$, and hence $$\max_{(x,\theta)\in \Gamma_{A_{k-1}}}\theta \ge \max_{(x,\theta)\in \Gamma_{K}}\theta.$$
	From \eqref{T1P1}, the inequality above and the definition of $\text{LB}_{k-1}$, 
	the following equations hold $$\max_{(x,\theta)\in \Gamma_{A_{k-1}}}\theta=\max_{(x,\theta)\in \Gamma_{K}}\theta=  \max_{x\in K}f(x) = \text{LB}_{k-1} = f(x^l),$$ for some $l\in \{0,1,\ldots,k-1\}$, and hence $x^l$ is optimal for \eqref{prob:pdsp}.
	
	Now, we show that Algorithm~\ref{Al-1} converges after a finite number of steps. Suppose $x^{k_1} = x^{k_2}$ for some $k_2 > k_1\ge 0$. Let $(x^{k_2},\theta^{k_2})\in \Gamma_{A_{k_2}}$, and $\theta^{k_2} > \text{LB}_{k_2}$. Then,
	$$
	\text{LB}_{k_2}\ge f(x^{k_1})=h(x^{k_2},x^{k_1}) \ge \theta^{k_2},
	$$
	which is a contradiction.
	
	This shows that Algorithm~\ref{Al-1} will not revisit a previous point. Since the set $K$ is finite, we must have finite convergence.
	\qed
\end{proof}

We now study the efficiency of the cutting planes by answering the question, at each iteration $k$, how many non-optimal solutions are eliminated by the cut $h(x,x^k) \le \theta$? We first establish the following elementary result.
\begin{lemma}
	\label{L1}
	Let $x^k\in K$ ($k\ge0$) be the iterate generated by Algorithm~\ref{Al-1} during iteration $k$. Then, for the subsequent iterations, we have
	\begin{equation}\label{GC1}
		0<	\theta^l - f(x^k)\le \ang{\nabla f(x^{k}), x^l -x^k},\quad \forall l> k.
	\end{equation}
	Consequently, if $\nabla f(x^{k}) \neq 0$, then for any $l > k$, it holds 
	\begin{equation}\label{GC2}
		\norm{x^t -x^k}\ge \frac{\text{\rm LB}_l- f(x^k)}{\norm{\nabla f(x^{k})}}, \quad \forall t \ge l.
	\end{equation}
\end{lemma}

\begin{proof}
	Consider $l>k$, and $(x^l,\theta^l)\in \Gamma_{A_{l}}$ such that $\theta^l > \text{LB}_l\ge f(x^k)$. Then, 
	$$
	\ang{\nabla f(x^k),x^l-x^k} +f(x^k) \ge \theta^l > \text{LB}_l\ge f(x^k). 
	$$
	The second assertion follows from \eqref{GC1} and the fact that for any $t\ge l$, $\text{LB}_t\ge \text{LB}_l$, and
	$$\norm{x^t -x^k}\ge \frac{\text{\rm LB}_t- f(x^k)}{\norm{\nabla f(x^{k})}} \ge \frac{\text{\rm LB}_l- f(x^k)}{\norm{\nabla f(x^{k})}}.$$
	\qed
\end{proof}

\begin{theorem}\label{element}
	At step $k$ of Algorithm~\ref{Al-1}, if ${\nabla f(x^k)} \neq 0$, and for some subsequent iterations $l>k$, there is $N_l\in \{0,1,\ldots,\min\{p,n-p\}\}$ such that
	\begin{align}\label{T5.1}
		\tfrac{\text{\rm LB}_l - f(x^k)}{\norm{\nabla f(x^k)}} > \sqrt{2N_l},
	\end{align}
	then from step $l$ onward the cutting plane $\theta \le h(x,x^k)$, removes at least 
	$
	\sum_{q=0}^{N_l}\binom{p}{q}\binom{n-p}{q}
	$
	binary points from the feasible set $K$, where $\binom{a}{b} = \frac{b!(a-b)!}{a!}$ for all $a,b\in \N$, and $a\ge b$. Consequently, if $N_l = \min\{p,n-p\}$, the cut $\theta \le h(x,x^k)$ removes all non-optimal solutions in $K$.
\end{theorem}

\begin{proof}
	Suppose \eqref{T5.1} holds for some $N_l\in \{0,1,\ldots,\min\{p,n-p\}\}$. It follows from Lemma~\ref{L1} that $\theta \le h(x,x^k)$ removes all points in $x\in K$ such that
	\begin{equation}\label{T5P1}
		\norm{x-x^k}^2\le 2N_l.
	\end{equation}
	Consider two sets of indices
	$$
	I_1:=\{i: \; x^k_i = 1\},\quad I_2:= \{i: \; x^k_i = 0\}.
	$$
	Since $x^k\in K$, we have $\abs{I_1} = p$ and $\abs{I_2} = n-p$. 
	For any $x\in K$, we consider 
	$$
	I_1(x):=\{i: \;x_i = 1\},\quad I_2(x):= \{i: \; x_i = 0\}.
	$$
	For any point $x\in K$, we have $\abs{I_1\setminus (I_1\cap I_1(x))} = \abs{I_2\setminus (I_2\cap I_2(x))}$, i.e., the number of indices of $x$ that are not selected in $I_1$ equals the number of indices of $x$ that are selected in $I_2$, and 
	$$
	\norm{x-x^k}^2 = \abs{I_1\setminus (I_1\cap I_1(x))} +  \abs{I_2\setminus (I_2\cap I_2(x))} =  2\abs{I_1\setminus (I_1\cap I_1(x))},
	$$
	i.e., the squared distance between $x$ and $x^k$ is the total number of indices differences between $x$ and $x^k$.
	We now start counting number of feasible solutions $x$ in $K$ such that the cardinality of ${I_1\setminus I_1\cap I_1(x)}$ is precisely $0,1,\ldots,N_l$ respectively. 
	
	For any $q\in \{0,1,\ldots,N_l\}$, consider $x\in K$ such that $I_1\setminus (I_1\cap I_1(x))$ is a subset of $I_1$ with cardinality $q$, and so there are $\binom{p}{q}$ such subsets of $I_1$. Similarly, $I_2\setminus (I_2\cap I_2(x))$ is a subset of $I_2$ with cardinality $q$, and so there are $\binom{p}{q}$ such subsets of $I_2$. Altogether, there are $\binom{p}{q}\binom{n-p}{q}$ feasible solutions $x\in K$ such that 
	$$
	\norm{x-x^k}^2= 2q.
	$$
	Therefore, there are precisely $
	\sum_{q=0}^{N_l}\binom{p}{q}\binom{n-p}{q}
	$ feasible points in $K$ that satisfy \eqref{T5P1}. This proves the assertion.
	\qed 
\end{proof}
When $\nabla f(x^k) = 0$, the cutting plane $h(x,x^k)\ge \theta$ becomes 
\begin{equation}
    \label{R1}
    f(x^k) \ge \theta.
\end{equation}
Hence, $f(x^k) \ge \max_{(x,\theta)\in \Gamma_{A_{k+1}}}\theta \ge \max_{(x,\theta)\in \Gamma_{K}}\theta = \max_{x\in K}f(x)$, and
the candidate $x^k$ is an optimal solution of \eqref{prob:pdsp}. 
The constraint \eqref{R1} implies that there will be no candidate solution found in iteration $k+1$, and hence Algorithm~\ref{Al-1} must terminate.
 
To finish this section, we explain how Algorithm~\ref{Al-1} can be integrated into the branch and cut procedure. 
This is achieved in standard MIP solvers using the \emph{callback} functionality. Callbacks allow for certain processes or algorithms to be implemented alongside general branch and bound or branch and cut procedures. In the case of Algorithm \ref{Al-1}, we begin by solving $(\text{LP}_{A_0})$ with a single cuts generated by some feasible solution. Then, a callback is used to add the associated cutting plane whenever an incumbent solution is found. Algorithm~\ref{Al-1} allows the MIP solver to preserve the information from previous solves and therefore generates only one search tree, improving the computational performance of the algorithm.

\section{Numerical Results}\label{NE}

We now look at Algorithm \ref{Al-1}'s performance on a range of test instances. The cutting plane algorithm was implemented in \texttt{CPLEX} Version 20.1 using the callback functionality. As explained above, callbacks allow for cuts to be added to the model during the general solve procedure, thus generating only one branch and cut search tree. The program was compiled using Microsoft Visual Studio and ran on a computer with a 1.8GHz Intel Core i7 processor with 16GB of RAM, using a single thread.

We begin by comparing the performance of the cutting plane algorithm with linear reformulation \eqref{prob:f3}, which is then solved using \texttt{CPLEX}. This linear reformulation was shown in \cite{Marti2022} to be competitive among other exact solvers.  \cite{Parreno2021} recently proposed a new large test set, labeled \texttt{GKD-d} and available within the \texttt{MDPLIB 2.0} test library \citep{MDLLIB2}.  In their recent review paper, the authors of \cite{Marti2022} suggest that this test set should replace previously used sets as they, for the first time, consider very large instances (up to $n=2000$) and include the original coordinate locations in the test set, thus allowing the results to be analysed similarly to \cite{Parreno2021}.

Test set \texttt{GDK-d} contains 70 matrices that define the Euclidean distances between randomly generated points with two coordinates in the range 0 to 100.  There are 10 matrices for each value $n=25,50,100,250,500,1000,2000$, the test is run with $p = \ceil{n/10}$ and $p = 2\ceil{n/10}$, totalling 140 test instances.  The test set is publicly available at \texttt{https://www.uv.es/rmarti/paper/mdp.html}.

Table \ref{tab:testA} compares the performance of the cutting plane algorithm against linear reformulation \eqref{prob:f3}, which is then solved using \texttt{CPLEX}, for each value $n$ in the test set \texttt{GDK-d}, over a time limit of 10 and 100 seconds.  For each $n$ and time limit the 20 tests are solved. We report the number of instances where the solver could confirm optimality and the average optimality gap. Over both time limits, the cutting plane algorithm appears vastly superior. On the 10-second time limit, the cutting plane algorithm was able to confirm optimality in all but two of the very large instances. However, over the 100-second limit, even these difficult instances were solved to optimality.  This represents a significant improvement when compared to \texttt{CPLEX}.

The same test was used as a part of the comprehensive empirical study conducted in \cite{Marti2022}. This study tested several state-of-art heuristic and meta-heuristic solvers on the 140 instances in \texttt{GDK-d}. The solvers were given 10 and 600-second time limits. For each time limit, the number of problems out of 140 solved to optimality, and the average deviation from the best-known solution were quoted. Over the 10-second time limit, the best heuristic solved 109 instances to optimality. In comparison, over a 10-second time limit, the cutting plane algorithm solved 138 instances to proven optimality. Over the 600-second time limit, the best heuristics solved 138 out of 140 instances, whereas the cutting plane algorithm solved all 140 within only a 100-second time limit.

\begin{table}
	\begin{center}
		\caption{Number of instances (out of 20) solved to optimality and average optimality gap for linear reformulation \eqref{prob:f3} (solved by \texttt{CPLEX}) and the cutting plane algorithm over a 10 and 100 second time limit.}
	    \label{tab:testA} 
		\begin{tabular}{r l rrrr l rrrr }
			& \phantom{a} & \multicolumn{4}{c}{10 seconds} & \phantom{a} & \multicolumn{4}{c}{100 seconds} \\
			\cmidrule{3-6} \cmidrule{8-11}
			&& \multicolumn{2}{c}{\texttt{CPLEX}} & \multicolumn{2}{c}{\texttt{Cutting Plane}} && \multicolumn{2}{c}{\texttt{CPLEX}} & \multicolumn{2}{c}{\texttt{Cutting Plane}} \\
			$n$ &&  \# opt &   \% gap &  \# opt & \% gap &&  \# opt &   \% gap &  \# opt & \% gap \\
			\midrule
			25 &&   20 &   0.00 &   20 & 0.00 &&  20 &   0.00 &   20 & 0.00 \\
			50 &&    0 &  47.27 &   20 & 0.00 &&  10 &  15.51 &   20 & 0.00 \\
			100 &&    0 &  90.48 &   20 & 0.00 &&   0 &  72.77 &   20 & 0.00 \\
			250 &&    0 & 130.12 &   20 & 0.00 &&   0 & 104.06 &   20 & 0.00 \\
			500 &&    0 & 162.95 &   20 & 0.00 &&   0 & 132.00 &   20 & 0.00 \\
			1000 &&    0 & 188.17 &   20 & 0.00 &&   0 & 159.03 &   20 & 0.00 \\
			2000 &&    0 &    -   &   18 & 0.01 &&   0 & 226.94 &   20 & 0.00 \\
			\bottomrule
		\end{tabular}  
	\end{center}
\end{table}

Table \ref{tab:testB} outlines the performances of the cutting plane algorithm on test set \texttt{GDK-d} in greater detail.  For each value $n$ and $p$, all ten instances were solved to optimality, and the average number of cuts added and solve time were reported. These results demonstrate the strength of the cuts as the number of cuts added remains very stable as the instance size grows. For instance, while 20.7 cuts were required on average for $n=25,p=3$, this only increased to 32.2 for $n=2000,p=400$.  The stability in the number of cuts required as the instance size grows is a testament to the cuts' strength and is a direct contributor to the fast runtimes. As the number of cuts remains low, the algorithm is able to outperform \texttt{CPLEX} comprehensively.

It is worth noting that the performance of the cutting plane algorithm seems to contradict a previously held notion about the difficulty of dispersion problems. \cite{Marti2022} state that as $p$ approaches $n/2$, a problem instance becomes harder due to the larger number of feasible solutions. While this may be true for many existing exact and heuristic solvers, this result was not observed for the cutting plane algorithm. As outlined before, every test instance in \texttt{GDK-d} is solved with $p = \ceil{n/10}$ and $p = 2\ceil{n/10}$.  The results in Table \ref{tab:testB} show that when $p$ is chosen as the larger value, fewer cuts are required on average, and therefore the problem is solved faster. This contradicts the statement in \cite{Marti2022} and shows that for the cutting plane algorithm, the run time does not increase as $p$ approaches $n/2$.

\begin{table}
	\begin{center}
		\caption{Performance of the cutting plane algorithm on the \texttt{GDK-d} test set.  For each $n$ and $p$, ten instances are solved to optimality and the average number of cuts and solve time in seconds are reported.}
	    \label{tab:testB}
		\begin{tabular}{rr l rrr lrrr}
			%\toprule
			&     &\phantom{a}&  \multicolumn{3}{c}{Cuts} & \phantom{a}& \multicolumn{3}{c}{Time (sec)} \\
			\cmidrule{4-6} \cmidrule{8-10}
			$n$ &   $p$ &&       Min &      Max &        Ave &&     Min &     Max &        Ave \\
			\midrule
			25 &   3 &&       11 &        37 &      20.7 &&   0.0156 &    0.0469 &    0.0297 \\
			25 &   6 &&        9 &        25 &      17.6 &&   0.0156 &    0.0312 &    0.0188 \\
			50 &   5 &&       21 &        49 &      31.5 &&   0.0469 &    0.1406 &    0.0734 \\
			50 &  10 &&       16 &        64 &      31.7 &&   0.0156 &    0.0938 &    0.0406 \\
			100 &  10 &&       19 &        76 &      41.3 &&   0.0625 &    0.4688 &    0.1781 \\
			100 &  20 &&       20 &        62 &      35.9 &&   0.0312 &    0.1875 &    0.0844 \\
			250 &  25 &&       28 &        93 &      60.2 &&   0.1875 &    1.9375 &    0.5984 \\
			250 &  50 &&       28 &        88 &      57.9 &&   0.1562 &    0.7031 &    0.3688 \\
			500 &  50 &&       22 &        70 &      51.8 &&   0.3750 &    1.6562 &    1.0391 \\
			500 & 100 &&       21 &       118 &      47.7 &&   0.3125 &    2.1094 &    0.8438 \\
			1000 & 100 &&       32 &        96 &      49.0 &&   2.0000 &    6.4062 &    3.1875 \\
			1000 & 200 &&       14 &        60 &      36.3 &&   0.9062 &    3.8281 &    2.2609 \\
			2000 & 200 &&       25 &        60 &      39.5 &&   6.0000 &   14.1406 &    9.5031 \\
			2000 & 400 &&       15 &        43 &      32.2 &&   3.5938 &   10.7031 &    7.8828 \\
			\bottomrule
		\end{tabular}
	\end{center}
\end{table}

It is important to note that \texttt{GDK-d} was generated using points with only two coordinates. We now explore the effect that increasing the number of coordinates has on the performance of the cutting plane algorithm. We do so by generating our test instances similar to \texttt{GDK-d}. As before, the matrix contains the Euclidean distance between randomly generated points with $s$ coordinates in the range 0 to 100, where $s$ is chosen in the range 2 to 20. For each $n=25,50,100$ we generate 10 instances with $p = \ceil{n/10}$.  Each instance is solved to optimality, and the run time and number of cuts required are reported.

\begin{figure}
	\begin{center}
		\input{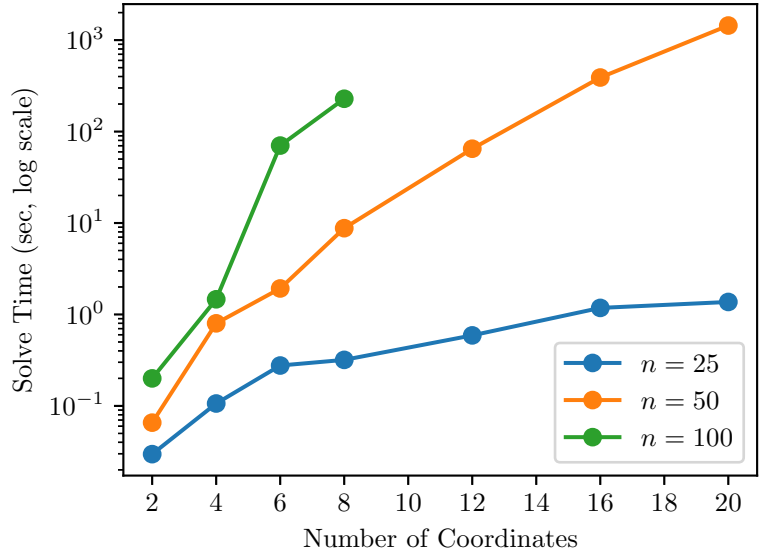}
	    \caption{Average solve time (on a logarithmic scale) for problems with increasing number of coordinates.}
	    \label{fig:testC_time}
	\end{center}
\end{figure}

Figure \ref{fig:testC_time} outlines the average solve time for each $n$ and $s$, where the solve time is plotted on a logarithmic scale. It clearly shows the exponential increase in solve time as the number of coordinates increases. For $s=2$, the solve time for $n=25,50,100$ remains below 1 second, which aligns with the results seen in Table \ref{tab:testB}.  However, as $s$ begins the increase, so does the solve time. For $n=50$, the solve time increases to an average of 1000 seconds for $s=20$.  The figure clearly shows the exponential increase in run time. For $n=100$, once $s$ was set greater than 8, the average solve time was far beyond 1000 seconds, and therefore for this problem size, we did not consider larger values of $s$. It is important to reinforce here that the problem dimensions are not increasing, as $n$ and $p$ are kept constant. Only the number of coordinates for the original points is increasing.

\begin{figure}
	\begin{center}
		\input{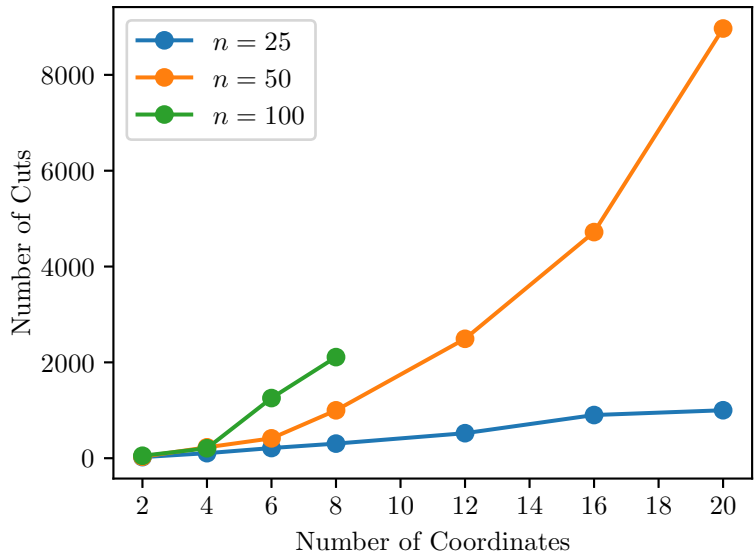}
	    \caption{Average number of cuts for problems with increasing number of coordinates.}
	    \label{fig:testC_cuts}
	\end{center}
\end{figure}

To understand the reason behind the drastic deterioration in performance when the number of coordinates increases, we should consider the number of cuts required to solve each instance. Figure \ref{fig:testC_cuts} shows the average number of cuts required to solve the 10 instances for each value $n,p$ and $s$.  Whereas Table \ref{tab:testB} demonstrated the stability in the number of cuts required for increasing $n$ when $s=2$, Figure \ref{fig:testC_cuts} shows how the number of cuts required increases exponentially as $s$ increases.  This represents a deterioration in the strength of each cut as $s$ increases. For example, for $n=50$, the number of cuts required increases to an average of over 8000 for $s=50$. This is far greater than the number of cuts required for even the very large test instances in \texttt{GDK-d} and implies that the cuts become incredibly weak for large values of $s$.  

\section{Conclusions and future work}
This paper presented a cutting plane algorithm for the $p$-dispersion-sum problem.   While the problem is inherently nonconvex, the cuts are shown to be appropriate, and the algorithm converges to the optimal solution. As the cuts can be applied directly to the original problem, the algorithm can avoid the unnecessary convexification steps commonly taken by integer quadratic solvers such as \texttt{CPLEX}. The algorithm's performance was evaluated on a recently published test library, where it was found to be vastly superior for Euclidean problems with a low number of coordinates. The cuts are tight for these location problems, allowing the algorithm to solve large instances quickly. However, as the number of coordinates grows, the cuts rapidly become inefficient, and the algorithm becomes far less effective. The reason for this is unclear, and future research should explore why the cuts become so poor for a larger number of coordinates and if there are ways to combat this challenge.

The cuts are appropriate because they do not remove feasible solutions and are upper planes of the objective in the domain of feasible solutions. However, the cuts are no longer valid if either the integrality condition or constraint \eqref{cts:p} is relaxed. This way, the problem can be considered concave/convex when the domain is restricted to purely discrete feasible solutions. This proof of concavity/convexity contrasts much of the previously conducted research into cutting plane algorithms for mixed-integer nonlinear programming, where concavity/convexity is usually assured by showing concavity/convexity when the integrality condition, constraint set, or both, are relaxed. Future research should (1) redefine concavity/convexity notions on the discrete domains and (2) uncover other well-known integer nonlinear problems where the problem is concave/convex on the domain of feasible solutions, but nonconcave/nonconvex when the domain is relaxed.
%, as cutting plane methods may be appropriate to solve these problems.

\section{Acknowledgment}
The authors are supported by the Australian Research Council through the Centre for Transforming Maintenance through Data Science (grant number IC180100030).
\addcontentsline{toc}{section}{References}

\bibliographystyle{apacite}
\bibliography{biblio}
\end{document}